\documentclass[a4paper,12pt]{article}

\usepackage{amsthm}
\usepackage{amsmath}
\usepackage{amssymb}

\numberwithin{equation}{section}

\newtheorem{thm}{Theorem}[section]
\newtheorem{prop}[thm]{Proposition}
\newtheorem{lem}[thm]{Lemma}
\newtheorem{dfn}[thm]{Definition}
\newtheorem{cor}[thm]{Corollary}
\newtheorem{ex}[thm]{Example}

\def\R{\mathbb{R}}
\def\Ri{\mathbb{R}\cup\{+\infty\}}
\def\N{\mathbb{N}}
\def\eps{\varepsilon}

\def\dist{\mathrm{dist\,}}

\def\n{|\hspace{-1pt}|\hspace{-1pt}|}
\def\be{\begin{equation}}
\def\ee{\end{equation}}
\def\n{|\hspace{-1pt}|\hspace{-1pt}|}

\def\ba{\begin{array}}
\def\ea{\end{array}}

\DeclareMathOperator{\diam}{diam}

\begin{document}

\title{Generic continuity of the perturbed minima of certain parametric optimization problems\footnote{The study  is supported by  the European Union-NextGenerationEU, through the National Recovery and Resilience Plan of the Republic of Bulgaria,  project  SUMMIT BG-RRP-2.004-0008-C01.}}
\date{}
\author{Hristina Topalova\thanks{Faculty of Mathematics and Informatics, Sofia University,   5, James Bourchier Blvd, 1164 Sofia, Bulgaria, e-mail: htopalova@fmi.uni-sofia.bg}\ \ and Nadia Zlateva\thanks{Faculty of Mathematics and Informatics, Sofia University,   5, James Bourchier Blvd, 1164 Sofia, Bulgaria, e-mail: zlateva@fmi.uni-sofia.bg}}
\maketitle
\begin{abstract}
	We show that in a quite general framework, the parameterized optimization problem can be so perturbed as to be generically well-posed.

    As an application, we provide a contribution to Ste\v{c}kin theory.


\end{abstract}

\section{Introduction}\label{sec:intro}

The article \cite{SimultOptim} by the same authors provides a new approach towards the classical, e.g. \cite{Danzig}, exploration of the behaviour of the  sets of minima of a parameterized minimization problem. This new approach is based on perturbation in the style of \cite{DGZ-article,DGZ}. The work \cite{SimultOptim}, focused on a sequence of minimization problems, hints that the method can perhaps be generalized; and the purpose of the present work is exactly this.

While it is not possible to have the desired result everywhere, see Example~\ref{ex:vime} below, in a generic sense it still holds.

We will now present our results in full details, and to this end we need quite a few technical definitions.

Our space framework is the metric spaces. Metric spaces are convenient, because -- unlike in normed spaces -- the subsets of metric spaces still fit the framework. If $(X,d)$ is a metric space then the closed ball with center $x\in X$ and radius $\varepsilon>0$ is denoted by $B_\varepsilon(x) := \{y\in X:\ d(y,x)\le\varepsilon\}$. For subsets $A,B\subset X$ we denote $d(A,B) := \inf\{d(x,y):\ x\in A,\ y\in B\}$. With a slight abuse of the notation we denote $d(a,B) := d(\{a\},B)$. For a nonempty bounded $A\subset X$,
$$
	\diam (A) := \sup\{d(x,y):\ x,y\in A\}.
$$
Our optimization framework is:
$$
	 \mathrm{Minimise}\ f_p(x)\ \mathrm{on}\ x\in X\ \mathrm{for}\ \mathrm{fixed}\ p\in P,
$$
where $(X,d)$ and $(P,\mu)$ are complete metric spaces, $p$ is regarded as parameter, and for each $p\in P$ the function $f_p:X\to\Ri$ is proper, lower semicontinuous and bounded below.

For a proper and bounded below function $f:X\to\Ri$ we denote the $\varepsilon$-argmin set by
$$
	\Omega_f(\varepsilon) := \{x\in X:\ f(x)\le \inf f + \varepsilon\}.
$$
The general problem under which the current work falls, is to investigate the continuity properties of the multivalued map
$$
	p \to \Omega_{f_p}(\varepsilon).
$$
Of course, to expect any continuity there ought to be a relationship between $f_p$'s for different $p$'s. What works in our case is what we will be calling \textit{uniform epi-continuity}. The term stems from the fact that for sequences it is an uniformized version of the epi-convergence, e.g. \cite[p.169]{beer}, or \cite[Chapter~7]{rocwets}. So, if $(X,d)$ and $(P,\mu)$ are complete metric spaces, the family of the functions $(f_p)_{p\in P}$ is called uniformly epi-continuous on $p\in P$, provided it satisfies the conditions \eqref{eq:cond-1} and \eqref{eq:cond-2} below:
\begin{eqnarray}
	\label{eq:cond-1}
	&&\forall p\in P,\ \forall x\in X,\ \forall \varepsilon>0\ \exists \delta>0 \mbox{ s.t.}\notag\\
	&&\forall q\in B_\delta(p)\ \exists  x_q \in B_\varepsilon(x) \mbox{ s.t.}\\	
	&&f_q(x_q) \le f_p(x) + \varepsilon,\notag
\end{eqnarray}
\begin{eqnarray}
	\label{eq:cond-2}
	&&\forall p\in P,\ \forall \varepsilon>0\ \exists \delta>0 \mbox{ s.t.}\notag\\
	&&f_q(x) \ge (f_p)_\varepsilon(x) - \varepsilon,\quad\forall x\in X,\ \forall q\in B_\delta(p).
\end{eqnarray}
Here, the regularization $(v)_\varepsilon$, where $\varepsilon>0$, of a bounded below function $v:X\to\Ri$ is
\begin{equation}
	\label{eq:v-eps}
	(v)_\varepsilon(x):=\inf \{ v(y): y\in B_\varepsilon(x)\}.
\end{equation}
Next, we axiomatically define feasible sets of functions that we will use for perturbation. We are not aiming for the most general possible axioms. Rather, we strive for a balance between generality and practicality. For a function $g:P\times X\to \R$ we write for consistency for a fixed $p\in P$
$$
	g_p(x) := g(p,x),\quad\forall x\in X,
$$
and then we can talk about a family of functions $(g_p)_{p\in P}$.

Denote by $BUC(X)$ the space of all bounded and uniformly continuous real valued functions on $X$. Equipped with the norm
$$
	\|v\|_\infty := \sup_{x\in X} |v(x)|,\quad\forall v\in BUC(X),
$$
$(BUC(X),\|\cdot\|_\infty)$  is a Banach space.
\begin{dfn}\label{pert_sp}
	Let $(X,d)$ and $(P,\mu)$ be  complete metric spaces. Let $(f_p)_{p\in P}$ be a family of proper, lower semicontinuous and bounded below functions from $X$ to $\Ri$. Let $f:P\times X\to \Ri$ be defined as $f(p,x):=f_p(x)$.
	
	A complete metric space
	$({\cal G} ,\rho)$ of functions $g:P\times X\to \R$ is called a \emph{perturbation space} on $P\times X$ for $f$ if the metric $\rho$ is translation invariant, that is $\rho(g_1+h,g_2+h)=\rho(g_1,g_2)$, and $\cal G$ is a convex cone with respect to the standard functional addition and scalar multiplication, and, moreover:
	\begin{itemize}
	\item[\emph{(i)}] $g_p\in BUC(X)$ for all $p\in P$;
	\item[\emph{(ii)}] the mapping
	\[
		p\to g_p
	\]
	from $P$ to $BUC(X)$, is continuous;
	\item[\emph{(iii)}]
	 $\forall p\in P$ there exists a constant $c(p)>0$ such that
	\begin{equation*}
	  \label{eq:dominate-uniform}
	  \|g_p\|_\infty \le c(p) \rho(g,0),\quad\forall g\in  {\cal G};
	\end{equation*}
	\item[\emph{(iv)}] $\forall p\in P$,   $\forall \varepsilon >0$ and $\forall g\in {\cal G} $, there exist $\delta>0$  and $g'\in {\cal G} $  such that
	\begin{equation*}
	  \label{eq:def:perturb:c}
	  \rho(g,g')<\varepsilon,  \text{ and } \diam(\Omega_{f_p+g'_p}(\delta)) < \varepsilon.
	\end{equation*}
	\end{itemize}
\end{dfn}

Now we can state our main result.
\begin{thm}
	\label{thm:main}
	Let $(P,\mu)$ and $(X,d)$ be complete metric spaces.  Let also $P$ be separable.
	
	Let $(f_p)_{p\in P}$ be a family of proper, lower semicontinuous and bounded below functions from $X$ to $\Ri$ which are uniformly epi-continuous, that is, they satisfy conditions \eqref{eq:cond-1} and \eqref{eq:cond-2}.

	Let $({\cal G},\rho)$ be a perturbation space on $P\times X$ for $f:P\times X\to \Ri$ which is defined as $f(p,x):=f_p(x)$.
	
	Let $W$ be a countable dense subset of $P$.
	
	Then there exists a dense $G_{\delta}$ subset $G$ of $\cal G$ such that for each $g\in G$ there exists a dense $G_{\delta}$ subset $P_g$ of $P$ such that $W\subset P_g$ and the following holds:

	For any fixed $p\in P_g$ the function $f_p+g_p$ attains its strong minimum at $x_{p,g}\in X$ and, moreover,
	\begin{equation}
		\label{eq:main-cont}
\forall p\in P_g,\ \forall \eps>0,\ \exists \delta>0 \ \colon \ 	\Omega_{f_q+g_q}(\delta) \subset B_\varepsilon(x_{p,g}),\quad \forall q\in B_\delta(p).
	\end{equation}
	Also, the functions
	$$
		p \to f_p(x_{p,g})\mbox{ and }p \to g_p(x_{p,g})
	$$
	are continuous on $P_g$.
\end{thm}

Recall that $\bar x$ is the strong minimum of $f:X\to\Ri$ if each minimizing sequence $(x_n)_{n\in\N}$, that is, such that $f(x_n)\to \inf_X f$, satisfies $x_n\to\bar x$. It is clear from Cantor~Lemma that if $X$ is a complete metric space and $f:X\to \Ri $ is proper and lower semicontinuous function, then
\begin{equation}
	\label{eq:str-min-iff}
	f\mbox{ attains strong minimum}\iff \lim_{\varepsilon\to0}\Omega_f(\varepsilon) = 0.
\end{equation}

Using $BUC(X)$ as a perturbation space, we get the following.
\begin{cor}
	\label{cor:select}
	Let $(P,\mu)$ and $(X,d)$ be complete metric spaces.  Let also $P$ be separable.
	
	Let $(f_p)_{p\in P}$ be a family of proper, lower semicontinuous and bounded below functions from $X$ to $\Ri$ which are uniformly epi-continuous, that is, they satisfy \eqref{eq:cond-1} and \eqref{eq:cond-2}.
	
	Let $W$ be a countable dense subset of $P$.

	Then for each $\varepsilon > 0$ there are a dense $G_{\delta}$ subset $V$  of $P$ such that $W\subset V$, and a continuous mapping
	$$
		\varphi : V \to X
	$$
	such that
	$$
		\varphi (p) \in \Omega_{f_p}(\varepsilon),\quad\forall p\in V.
	$$
\end{cor}
The next simple example shows that even in the case when $P$ an interval, one cannot expect that the conclusion of Corollary~\ref{cor:select} -- and, therefore, that of Theorem~\ref{thm:main} -- would always hold for all $p\in P$. Of course, that could be the case under additional assumption on $f_p$'s like convexity, see e.g. \cite{dev-pro,pando,pando2,vesely}, but we do not cover these results here.

\begin{ex}	\label{ex:vime}
	Let $P=X=[0,1]$ and let $f_p(x) := (1-p)(3x - 1)$ for $x\in[0,1/3]$,  $f_p(x) := 0$ for $x\in[1/3,2/3]$, and $f_p(x) := p(2 - 3x)$ for $x\in[2/3,1]$. Then for each $\varepsilon\in(0,1/2)$ we have that $\Omega_{f_0}(\varepsilon)\subset[0,1/3]$, $\Omega_{f_1}(\varepsilon)\subset[2/3,1]$, and $\Omega_{f_p}(\varepsilon)\cap(1/3,2/3)=\varnothing$ for all $p\in[0,1]$. Therefore, there cannot be a continuous selection of the map $p\to\Omega_{f_p}(\varepsilon)$.
\end{ex}
Of course, if $P$ is countable, then we can take $W=P$ and get the conclusion of Theorem~\ref{thm:main} on the whole $P$. For example, let $P=\N\cup\{\infty\}$ with $\mu(i,j) := |1/i-1/j|$ and $\mu(i,\infty) := 1/i$. In this case continuity of the strong minima means convergence of the sequence $(x_n)_{n\in\N}$ to $x_\infty$ and we fully recover Theorem~1.2 of \cite{SimultOptim} as a corollary to Theorem~\ref{thm:main}.

As a non-trivial application of our result we elaborate on Ste\v{c}kin theory, e.g. \cite{Stechkin,EfimovStechkin} (see also \cite{Cobzas} and \cite[pp.1223]{Tbook} for a modern and comprehensive review).

Let $(E,\|\cdot\|)$ be a Banach space. Denote by $B_E$ and $S_E$ the unit ball and the unit sphere of $E$, respectively. Let $M$ be a nonempty closed, convex and bounded subset of $E$. Obviously, with the canonical metric $\mu(x,y):= \|x-y\|$  both $(E,\mu)$ and $(M,\mu)$ are complete metric spaces.

Consider the set ${\cal N}$ of all bounded seminorms on $E$, i.e. the functions $\nu:E\to \R$ such that
\begin{itemize}
\item [(j)]$\nu(tx)=|t|\nu(x)$, $\forall x\in E$;
\item [(jj)]$\nu(x+y)\le \nu(x)+\nu (y)$, $\forall x,y\in E$;
\item [(jjj)]$k_\nu = \sup_{x\in B_E} |\nu(x)|$ is finite, so $\nu(x)\le k_\nu\| x\|$, $\forall x\in E$.
\end{itemize}
Equipped with the metric
\begin{equation}\label{rho}
\rho(\nu_1,\nu_2):=\sup_{x\in B_E} |\nu_1(x)-\nu_2(x)|,
 \end{equation}
$({\cal N},\rho)$ is a complete metric space. The set ${\cal N}_0$ of all equivalent norms on $E$ is an open dense set in ${\cal N}$, see Lemma~\ref{lem:N0-open}.

\begin{prop}
	\label{pro:stech-ps}
Let $E$ be a Banach space.	The cone $({\cal G},\rho)$ of the functions
	$$
		{\cal G}:=\{ g:E \times E\to \R \mbox{ s.t. }g(p,x)=\nu (p-x),\ \nu \in {\cal N}\}
	$$
	equipped with the metric
	$$
		\rho(g_1,g_2):=\rho(\nu_1,\nu_2),\text{ where }g_i(p,x)=\nu_i(p-x),\ i=1,2,
	$$
	  is a perturbation space on $P\times M$ for the function $f(p,x)\equiv 0$ for all $p\in E$ and $x\in M$.
\end{prop}
Note that here the functions $(f_p)_{p\in E}$ do not actually depend  on $p$ (and neither on $x$ for that matter). However, the functions in the perturbation space do depend on $p$. It is because of this  example and the next result that we added the additional complexity (in \cite{SimultOptim} the perturbation functions do not depend on the parameter $p$).

A subtlety that arises here is that different seminorms may have identical restrictions on the set $M$, but this  does not cause a problem. First of all, Definition~\ref{pert_sp} allows for $g_1,g_2\in\cal G$ such that $\rho(g_1,g_2) > 0$, but $g_1(p,x) = g_2(p,x)$ for all $p\in P$ and $x\in X$. Second, in our concrete example, if $\nu_1(p-x) = \nu_2(p-x)$ for all $p\in E$ and $x\in M$, then obviously $\nu_1\equiv\nu_2$.

The Ste\v{c}kin-like result we get  is as follows.
\begin{thm}
	\label{thm:stech}
	Let $(E,\|\cdot\|)$ be a separable Banach space, and let $M\subset E$ be a nonempty closed convex and bounded set in $E$. For any dense countable set $W\subset E$ and
	any equivalent norm $|\cdot|$ on $E$,  there exist an arbitrary close to $|\cdot|$
	equivalent norm $\n\cdot\n$ and a dense $G_\delta$ set $V\supset W$ such that the problem
\begin{equation}\label{mp}
		\min _{x\in M}\n p-x\n
\end{equation}
	is well-posed for any   $p\in V$. In particular, the metric projection on $M$ with respect to the norm $\n\cdot\n$  is single valued and continuous on $V$.
\end{thm}
Recall that for $Z\subset E$ and $f:Z\to \Ri$ the minimization problem $\min _Z f$ is \emph{well-posed} if the restriction of $f$ to $Z$   attaining strong minimum in $Z$. Note that this result is original. In the classical Ste\v{c}kin theory the norm is fixed, but it has some nice properties. Obviously, if the norm is uniformly convex, then the metric projection on the closed convex set $M$ will be well-posed everywhere. But not every separable Banach space admits an uniformly convex renorm. \emph{Locally uniformly convex} renorms are available, but in general they do not guarantee well-posedness of the minimization problem \eqref{mp} anywhere, see \cite[p. 276]{Cobzas}. In \cite[Theorem~4.11]{Pando-st} it is proved that for any closed set and any equivalent norm there is an equivalent norm close to the latter, such that the metric projection is generically well-posed. The new feature here is that we can fix in advance a dense countable set $W$ and then get the renorm in such a way that the metric projection will be well-posed and, therefore, continuous on $W$ (in fact it will be continuous on a dense $G_\delta$ set containing $W$).

The work is organized as follows. In Section~\ref{sec:perturb} we adapt the perturbation method from \cite{SimultOptim} to our slightly more general Definition~\ref{pert_sp}. In Section~\ref{sec:epic} we show that the uniform epi-continuity implies that the set of parameters where strong minimum is attained, is $G_\delta$. These two combined easily yield the proof of our main Theorem~\ref{thm:main} in Section~\ref{sec:mpr}. The final Section~\ref{sec:stech} is devoted to our contribution to Ste\v{c}kin theory. 

\section{Perturbation method}\label{sec:perturb}

In \cite{DGZ-article,DGZ} the authors established not only a variant of Borwein-Priess Smooth Variational Principle~\cite{B-P}, but also gave a very general  scheme for creating  variational principles or, what some prefer to say, perturbation methods. This scheme has been followed e.g. in \cite{dev-pro,iv-zla-jca,iv-zla-jota,pertOrl} and we will follow it here as well.

In this section, we demonstrate the application of the perturbation method for a fixed parameter $p\in P$. Specifically, we present the following result, which is a minor adaptation of \cite[Theorem~3.1]{SimultOptim}. Due to our slightly more general definition of a perturbation space, we provide a complete proof here.

\begin{thm}[\textbf{Variational Principle}]\label{thm:VP1}
    Let $(P,\mu)$ and $(X,d)$ be  complete metric spaces.

    Let the functions  $f_p: X\to \Ri$ be  proper, lower semicontinuous and bounded below for all $p\in P$. Let $({\cal G},\rho)$ be a perturbation space for the function $f:P\times X\to \Ri$   defined as $f(p,x):=f_p(x)$.

    Let $p\in P$ be fixed. Then there exists a dense $G_\delta$ subset $U$ of ${\cal G}$ such that for each $g\in U$ the function $f_p+g_p$  attains strong minimum on $X$.
\end{thm}
For the proof we need the following useful Lemma, cf. \cite[Lemma~4]{iv-zla-jca}, which reveals certain continuity-like property of the $\varepsilon$-argmin map.

\begin{lem}\label{lem:cont-eps}
    Let $(X,d)$ be a metric space and let $f:X\to\Ri$ be proper and bounded below function. Let $\varepsilon > 0$.
    If the bounded function $g:X\to\R$ is such that
    $$
        \|g\|_\infty < \varepsilon/3,
    $$
    then
    $$
        \Omega_{f+g}(\varepsilon/3) \subset \Omega_f(\varepsilon).
    $$
\end{lem}
\begin{proof}
    Let $\delta := \varepsilon/3$, so $\|g\|_\infty < \delta$. Let $x\in \Omega_{f+g}(\delta)$.

    It is clear that
    $$
        \inf (f+g) < \inf f + \delta,
    $$
    so, $f(x) < f(x) + g(x) + \delta \le \inf(f+g) + 2\delta < \inf f + 3\delta$. That is, $x\in \Omega_f(\varepsilon)$.
\end{proof}
\begin{proof}[\textbf{\emph{Proof of Theorem~\ref{thm:VP1}}}]
    Fix $p\in P$ and consider for $n\in\N$ the subset $M_{n}$ of ${\cal G}$ defined by
    \begin{equation}\label{eq:an-def}
        M_{n}:=\left\{ g\in {\cal G}:  \exists\ t>0  :\ \diam\left(\Omega_{f_p+g_p}(t)\right) < \frac{1}{n}\right\}.
    \end{equation}
    We will show that $M_{n}$ is dense and open in $({\cal G},\rho)$. Then by Baire Category Theorem
    $$
	    U := \bigcap_{n\in \N} M_{n}
    $$
    will be dense in ${\cal G}$ (it is $G_\delta$ by definition). For each $g\in U$ the definition of $M_{n}$ shows that $\lim_{\eps \to 0} \diam (\Omega_{f_p+g_p}(\eps)) =0$, hence $f_p+g_p$ attains  strong minimum, see \eqref{eq:str-min-iff}.

    Fix arbitrary $n\in \N$.

    To show that $M_n$ is indeed dense, take arbitrary $h\in {\cal G}$ and  $\varepsilon \in (0,1/n)$.

    Let $\delta > 0$ and $g'\in {\cal G}$ be provided by (iv) of Definition~\ref{pert_sp} for $h$. So, $\rho (h,g')<\eps$ and $\diam (\Omega_{f_p+g_p'}(\delta))<\eps<\dfrac{1}{n}$.
    This means that $g'\in M_{n}$. Since $\rho(h,g')<\varepsilon$, the distance from $h$ to $M_{n}$ is smaller than $\eps$. In other words, $M_{n}$ is dense in~${\cal G}$.

    To show that $M_n$ is open, take an arbitrary $g\in M_{n}$. This means by \eqref{eq:an-def} that there is $\varepsilon > 0$ such that
    $$
        \diam\left(\Omega_{f_p+g_p}(\varepsilon)\right) < \frac{1}{n}.
    $$
    Let $\delta > 0$ be such that
    $$
        \delta < \frac{\varepsilon}{3c(p)},
    $$
    where $c(p)$ is the constant from  Definition~\ref{pert_sp}(iii). Let $g'\in B_\delta(g)$ be arbitrary. Consider
    $$
        h := g_p' - g_p.
    $$
    By   Definition~\ref{pert_sp}(iii) it follows that $\|h\|_\infty \le c(p) \rho(g_p',g_p) \le c(p)\delta < \varepsilon/3$. Thus Lemma~\ref{lem:cont-eps} gives $\Omega_{f_p+g_p + h}(\varepsilon/3) \subset \Omega_{f_p+g_p}(\varepsilon)$. But $f_p+g_p + h = f_p+g_p'$, so
    $$
        \diam\left(\Omega_{f_p+g_p'}(\varepsilon/3)\right)\le \diam\left(\Omega_{f_p+g_p}(\varepsilon )\right) <   \frac{1}{n} \Rightarrow g'\in M_n.
    $$
\end{proof} 

\section{Epi-continuity}\label{sec:epic}

    From Theorem~\ref{thm:VP1} it is clear how to get a countable set of parameters such that strong minimum to be attained. Here we will show that the uniform epi-continuity guarantees that the set to be $G_\delta$.

    \begin{prop}\label{prop:epic}
        Let $(P,\mu)$ and $(X,d)$ be  complete metric spaces. Let the functions  $h_p: X\to \Ri$ be  proper, lower semicontinuous and bounded below for all $p\in P$. Let also $(h_p)_{p\in P}$ be uniformly epi-continuous, that is, they satisfy \eqref{eq:cond-1} and \eqref{eq:cond-2}.

        Let $Q\subset P$ be the (possibly empty) set of these $p\in P$ for which $h_p$ attains strong minimum, say $x_p\in X$.

        Then $Q$ is a $G_\delta$ set and, moreover,
        \begin{equation}
            \label{eq:epic-prop-concl}
            \forall p\in Q,\ \forall \varepsilon > 0,\  \exists \delta >0\colon \
            \Omega_{h_q}(\delta) \subset B_{\varepsilon}(x_p),\quad \forall q\in B_\delta(p).
        \end{equation}
    \end{prop}

    We start with the following result somehow similar to Lemma~\ref{lem:cont-eps}.
    \begin{lem}\label{lem:epic}

        Let $(P,\mu)$ and $(X,d)$ be  complete metric spaces. Let the functions  $h_p: X\to \Ri$ be  proper, lower semicontinuous and bounded below for all $p\in P$. Let also $(h_p)_{p\in P}$ be uniformly epi-continuous, that is, they satisfy the conditions \eqref{eq:cond-1} and \eqref{eq:cond-2}.

        If for some $p\in P$, $\varepsilon > 0$ and $r > 0$,
        $$
            \diam\Omega_{h_p}(\varepsilon) < r,
        $$
        then there exists $\delta>0$ such that
        $$
            \diam\Omega_{h_q}(\delta) < 5r,\quad\forall q\in B_\delta(p).
        $$
    \end{lem}
    \begin{proof}
        Fix $p\in P$, $\varepsilon > 0$ and $r > 0$ such that $\diam\Omega_{h_p}(\varepsilon) < r$, and let $\mu>0$ be such that
        $$
            \mu \le \min \{r,\varepsilon/4\}.
        $$
        Let us fix $x\in \Omega_{h_p}(\mu)$. Let $\delta \in (0,\mu)$ be such that \eqref{eq:cond-1} and \eqref{eq:cond-2} hold for the fixed $p$,  $x $,  and $ \mu$.  We claim that
        \begin{equation}
            \label{eq:epiclaim}
            \Omega_{h_q}(\delta) \subset B_{2r}(x),\quad\forall q\in B_\delta(p)
        \end{equation}
    which is enough to get the conclusion.

    Indeed, take arbitrary  $y\in X$  such that $d(y,x) > 2r$.

        Since $\mu <\eps$, $x\in \Omega_{h_p}(\varepsilon)$ and the diameter of the latter set is smaller than~$r$ by assumption, we have that $B_{r}(y)\cap \Omega_{h_p}(\varepsilon)=\varnothing$ and because of $\mu\le r$,
        $$
            (h_p)_\mu(y) \ge (h_p)_r(y) > \inf h_p +\varepsilon \ge h_p(x) - \mu + \varepsilon \ge h_p(x) + 3\mu.
        $$
        From \eqref{eq:cond-2} it follows that
        $$
            h_q(y) \ge (h_p)_\mu(y)-\mu > h_p(x) + 3\mu-\mu = h_p(x) + 2\mu,\quad\forall q\in B_\delta(p).
        $$
        From \eqref{eq:cond-1} on the other hand, we get for each $q\in B_\delta(p)$ a $x_q\in X$ such that $h_q(x_q) \le h_p(x) + \mu$. Combining with the above inequality we get  that for all $q\in B_\delta(p)$,
        $$
    h_q(y) - \inf h_q \ge        h_q(y) - h_q(x_q) \ge        h_q(y) - h_p(x)-\mu >2\mu -\mu=\mu.$$
    Therefore, $ y \not\in \Omega_{h_q}(\mu)$ and since $\delta <\mu$, $y \not\in \Omega_{h_q}(\delta)$. Thus  \eqref{eq:epiclaim} is verified, because $y\not\in B_{2r}(x)$ was arbitrary.
    \end{proof}

    \begin{proof}[\textbf{\emph{Proof of Proposition~\ref{prop:epic}}}]
        Let
        $$
            A_n := \{p\in P:\ \exists \delta >0:\ \diam \Omega_{h_p}(\delta) < 5^{-n}\}.
        $$
        From \eqref{eq:str-min-iff} it follows that
        $$
            Q = \bigcap_{n\in\N} A_n.
        $$
        From Lemma~\ref{lem:epic} it holds that for any $p\in A_n$ there is $\delta '>0$ such that $B_{\delta '}(p) \subset A_{n+1}$, so $p\in A_{n+1}^\circ$, where $A^\circ$ is the interior of the set $A\subset X$. Therefore,
        $$
            A_{n+1} \subset A_n^\circ\subset A_n,\quad\forall n\in\N,
        $$
        hence
        $$
            Q = \bigcap_{n\in\N} A_n^\circ,
        $$
        and then $Q$ is a $G_\delta$ set.

        To establish  \eqref{eq:epic-prop-concl} take $p\in Q$ (if any) and $\eps >0$. As $\diam \Omega_{f_p}(\delta')\to 0$ whenever $\delta '\searrow 0$, there is $\delta '>0$ such that $\diam \Omega_{f_p}(\delta')<\eps/5$. By Lemma~\ref{lem:epic} there is $\delta >0$ such that
        $\diam \Omega_{f_q}(\delta) <\eps$ and \eqref{eq:epic-prop-concl} holds.
    \end{proof}

\section{Proof of the main result}
    \label{sec:mpr}
 To prove our main Theorem~\ref{thm:main}, we will need the following simple proposition. It  says in essence that the uniform epi-continuity implies continuity of the value function. This is one of the reasons why the   concept of uniform epi-convergence worths consideration.

    \begin{prop}
        \label{pro:epi-cont-V}
        Let $(X,d)$ and $(P,\mu)$ be complete metric spaces. Let the family of proper and bounded below functions $(f_p)_{p\in P}$ be uniformly epi-continuous on $p\in P$. Then the value function
        $$
            V(p) := \inf_{x\in X} f_p(x),\quad\forall p\in P,
        $$
        is continuous.
    \end{prop}
    \begin{proof}
        Fix $p\in P$ and $\varepsilon > 0$.

        Let $x_p\in X$ be such that $f_p(x_p) < V(p) + \varepsilon/2$. From \eqref{eq:cond-1} it follows that there is $\delta_1 > 0$ such that for all $q\in P$ such that $0 < \mu(p,q) < \delta_1$ there is $x_q\in X$ such that $f_q(x_q) \le f_p(x_p) + \varepsilon/2$. This means that
        $$
            V(q)\le f_q(x_q)  \le f_p(x_p) + \varepsilon/2 \le V(p) + \varepsilon,\quad\forall q\in B_{\delta_1}(p).
        $$

        On the other hand, from \eqref{eq:cond-2} there is $\delta_2 > 0$ such that for all $q\in B_{\delta_2}(p)$ and all $x\in X$ we have $f_q(x) \ge (f_p)_{\varepsilon} - \varepsilon$. Thus,
        $$
            V(q) \ge \inf (f_p)_{\varepsilon} - \varepsilon \ge V(p) - \varepsilon,\quad\forall q\in B_{\delta_2}(p).
        $$

        So, $p\to V(p)$ is continuous.
    \end{proof}
    The fundamental equivalence between Heine and Cauchy definitions (of convergence, continuity, etc.) holds for uniform epi-continuity as well. More precisely, if $(X,d)$ and $(P,\mu)$ are complete metric spaces, then the family of functions $(f_p)_{p\in P}$ is uniformly epi-continuous on $p\in P$, that is, it satisfies the conditions \eqref{eq:cond-1} and \eqref{eq:cond-2}, if and only if, for each $p\in P$ and each sequence $(p_n)_{n\in\N}$ convergent to $p$ (that is, $\mu(p,p_n)\to0$, as $n\to\infty$) the sequence of functions $\phi_n := f_{p_n}$ converges to $\phi_\infty := f_{p}$ in the sense of \cite{SimultOptim}, that is, $(\phi_n)_{n\ge1}$ satisfy the conditions (3) and (4) of \cite{SimultOptim}. In this way we can just reformulate \cite[Proposition~2.2]{SimultOptim}. However, for reader's convenience we will present a proof. In essence, the next proposition establishes a condition under which the uniform epi-continuity is preserved under addition.
    \begin{prop}
        \label{pro:sum-epi}
        Let $(X,d)$ and $(P,\mu)$ be complete metric spaces. Let the family of proper and bounded below functions $(f_p)_{p\in P}$ be uniformly epi-continuous on $p\in P$.

        Let the family $(g_p)_{p\in P}\subset BUC(X)$ be such that $p\to g_p$ is continuous from $P$ to $BUC(X)$ (in other words, $(g_p)_{p\in P}$ satisfy (i) and (ii) of Definition~\ref{pert_sp}).

        Then the family $(f_p + g_p)_{p\in P}$ is uniformly epi-continuous on $p\in P$.
    \end{prop}
    \begin{proof}
        From the conditions on $(g_p)_{p\in P}$ it easily follows that the function $(p,x)\to g_p(x)$ is continuous on both variables with uniformity on $x$, that is,
        \begin{equation}
            \label{eq:g-cont}
  \begin{array}{l}
  \forall p\in P,\ \forall \varepsilon>0, \exists \delta>0 \colon  \\[8pt]
  |g_q(y)-g_p(x)|<\varepsilon,\quad\forall q\in B_\delta(p),\ \forall x\in X,\ \forall y\in B_\delta(x).
  \end{array}
        \end{equation}
        Indeed, fix $p\in P$ and $\eps >0$. Let $\delta>0$ be such that $\|g_q-g_p\|_\infty < \varepsilon/2$ for all $q\in B_\delta(p)$, and also $|g_p(y)-g_p(x)| < \varepsilon/2$ for all $x\in X$ and $y\in B_\delta(x)$. Then for all $x\in X$ and $(q,y)\in B_\delta(p)\times B_\delta(x)$, we have $|g_q(y)-g_p(x)| \le |g_q(y)-g_p(y)| + |g_p(y)-g_p(x)| < \varepsilon$.

        Now, to verify \eqref{eq:cond-1}, fix $p\in P$, $x\in X$ and $\varepsilon>0$. Let $\xi \in (0,\varepsilon/2)$ be so small that for all $q\in B_\xi(p)$ and all $y\in B_\xi(q)$ we have $|g_q(y) - g_p(x)| < \varepsilon/2$, see \eqref{eq:g-cont}.

        Because $(f_p)_{p\in P}$ satisfy \eqref{eq:cond-1}, we can find $\delta\in(0,\xi)$ so small that for any $q\in B_\delta(p)$ we can fix a $x_q\in B_\xi(x)$ such that $f_q(x_q) \le f_p(x) + \xi$.

        Then for all $q\in B_\delta(p)$ and $y\in B_\delta(x)$ we have $(f+g)_q(x_q) - (f+g)_p(x) = (f_q(x_q)-f_p(x)) + (g_q(x_q)-g_p(x)) < \xi + \varepsilon/2 < \varepsilon$, that is, $(f_p+g_p)_{p\in P}$ satisfy \eqref{eq:cond-1}.

        Turning now to \eqref{eq:cond-2}, fix $p\in P$ and $\varepsilon>0$. Because $g_p$ is uniformly continuous, there is $\xi\in(0,\varepsilon/3)$ such that $|g(y)-g(x)|<\varepsilon/3$ for all $x$ and $y$ in $X$ such that $d(y,x)<\xi$. Then it is easy to check that
        \begin{equation}
            \label{eq:inf-eps}
            (f_p+g_p)_\xi (x) \le (f_p)_\xi(x)+g_p(x) + \varepsilon/3,\quad\forall x\in X.
        \end{equation}
        Indeed, $(f_p+g_p)_\xi (x)=\inf\{f_p(y)+g_p(y):\ y\in B_\xi(x)\} \le \inf\{f_p(y)+g_p(x):\ y\in B_\xi(x)\} + \sup\{g_p(y)-g_p(x):\ y\in B_\xi(x)\}$, and the latter supremum does not exceed $\varepsilon/3$.

        Let $\delta\in(0,\xi)$ be so small that for all $q\in B_\delta(p)$ we have that $\|g_q-g_p\|_\infty < \xi$, and also
        $$
            f_q(x) \ge (f_p)_\xi(x) - \xi,\quad\forall x\in X.
        $$
        This is possible, because $(f_p)_{p\in P}$ satisfy \eqref{eq:cond-2}. Then for all $q\in B_\delta(p)$ and all $x\in X$ we have $f_q(x) + g_q(x) \ge (f_p)_\xi(x) + g_p(x) - 2\xi$. Taking into account \eqref{eq:inf-eps}, we get $f_q(x) + g_q(x) \ge (f_p+g_p)_\xi (x) - \varepsilon$. Because trivially $(f_p+g_p)_\xi \ge (f_p+g_p)_\varepsilon$, we conclude that $(f_p+g_p)_{p\in P}$ satisfy \eqref{eq:cond-2}.
    \end{proof}

    The proof of our main result is now straightforward.

    \begin{proof}[\textbf{\emph{Proof of Theorem~\ref{thm:main}}}]
        First, we apply Theorem~\ref{thm:VP1} for each fixed $p\in W$ to get a dense $G_\delta$ subset $U_p$ of $\cal G$ such that $W\subset U_p$ and for each $g\in U_p$ the function $f_p+g_p$ attains strong minimum on $X$.

        We will show that
        $$
            G := \bigcap_{p\in W} U_p
        $$
        satisfies the conclusions of Theorem~\ref{thm:main}.

        What is immediate from Baire Theorem, is that $G$ is a dense and $G_\delta$ subset of $\cal G$, because $U_p$'s are such and $W$ is countable.

        Fix a $g\in G$. Then for each $p\in W$ the function $f_p+g_p$ attains strong minimum on $X$.

        From Proposition~\ref{pro:sum-epi} we have that $(f_p+g_p)_{p\in P}$ is uniformly epi-continuous on $P$, so Proposition~\ref{prop:epic} applied for $h_p=(f+g)_p$ gives that the set $P_g$ of those $p\in P$ for which $(f+g)_p$ attains strong minimum on $X$, is a $G_\delta$ set. Since we already know that the dense set $W\subset P_g$, the  set $  P_g$ is a dense $G_\delta$ set.

        Let for $p\in P_g$ the point $x_{p,g}\in X$ be the strong minimum of $(f+g)_p$, i.e.
        $$
            (f+g)_p (x_{p,g}) = \min _{x\in X} (f+g)_p(x),\quad\forall p\in P_g.
        $$
        The conclusion \eqref{eq:epic-prop-concl} of Proposition~\ref{prop:epic} for $h_p=(f+g)_p$ is identical to  \eqref{eq:main-cont}, so  \eqref{eq:main-cont} holds. Hence,  we only have to check the continuity of the functions $p\to f_p(x_{p,g})$ and $p\to g_p(x_{p,g})$ on $P_g$. To this end, note that \eqref{eq:main-cont} implies that the mapping $p\to x_{p,g}$ is continuous on $P_g$. From \eqref{eq:g-cont}, which $g$ certainly satisfies, it follows that $p\to g_p(x_{p,g})$ is continuous on $P_g$. Because $p\to (f+g)_p(x_{p,g})$ is uniformly epi-continuous on $P$ we can use  Proposition~\ref{pro:epi-cont-V} to conclude.
    \end{proof}

    \begin{proof}[\textbf{\emph{Proof of Corollary~\ref{cor:select}}}]
        It is easy to see that the set of functions $BUC(X)$, which appears in this fashion already in the seminal \cite{DGZ-article}, see also \cite{DGZ}, is a perturbation space for $(f_p)_p\in P$.

      Indeed, since $g\in BUC(X)$ does not depend on $p$, we only have to check  Definition~\ref{pert_sp}(iv). To this end we consider the functions
        $$
            u_{a,\beta,\gamma}(x) := \begin{cases}
                (\beta/\gamma)d(x,a),& x\in B_\gamma(a),\\
                \beta,& x\not\in B_\gamma(a),
            \end{cases}
        $$
        where $a\in X$ and $\beta,\gamma>0$. Obviously, $u_{a,\beta,\gamma}\in BUC(X)$ and $\|u_{a,\beta,\gamma}\|_\infty = \beta$.

        Now, fix $p\in P$, $\eps >0$ and $g\in BUC(X)$. Let us agree to write $f$ instead of $f_p$ for short.
        Let us take $a \in \Omega_{f+g}(\varepsilon/4)$. For $\beta=\eps$ and $\gamma=\eps/2$ define
        $$
            g' := g + u_{a,\varepsilon,\varepsilon/2},
        $$
        to get  $g'\in BUC(X)$ such that $\|g-g'\|_\infty = \varepsilon$. If $d(x,a) > \eps/2$ for some $x\in X$, then
        $$(f+g')(x) = (f+g)(x) + u_{a,\varepsilon,\varepsilon/2}(x) = (f+g)(x) + \varepsilon \ge \inf_X (f+g) + \varepsilon.$$
        But
        $$ (f+g')(a) = (f+g)(a) \le \inf_X (f+g) + \varepsilon/4,$$
         which combined with the above inequality gives
        $$(f+g')(x)\ge  (f+g')(a)+3\eps/4 > \inf_X (f+g') +\eps/2.$$
        Thus, $x\not\in \Omega_{f+g'}(\varepsilon/2)$. So, $\Omega_{f+g'}(\varepsilon/2)\subset B_{\eps/2}(a)$, hence  $\diam \Omega_{f+g'}(\varepsilon/2) \le \eps$, and Definition~\ref{pert_sp}(iv) is fulfilled with $\delta=\eps/2$.

        Having checked that $BUC(X)$ is a perturbation space for $(f_p)_{p\in P}$, we can now apply Theorem~\ref{thm:main} to get for a fixed $\varepsilon > 0$, a function $g \in BUC(X)$ such that $\|g \|_\infty < \varepsilon/2$ and $f_p + g$ attains its strong minimum on $X$ for all $p$  in a dense $G_\delta$ subset $P_g$ of $P$. Let
        $\varphi(p)\in X$ be such that
        $$
            (f_p + g)(\varphi(p)) = \min_X (f_p + g),\quad \forall p\in  P_g.
        $$
        From \eqref{eq:main-cont} it follows that $p\to \varphi (p)$ is continuous on $P_g$. Since $\|g \|_\infty < \varepsilon/2$, it is immediate that $\varphi (p) \in \Omega_{f_p}(\varepsilon)$ for all $p\in P_g$. It is enough to take $V=P_g$ to conclude.
    \end{proof}

\section{On Ste\v{c}kin theory}\label{sec:stech}

Here we will show how Theorem~\ref{thm:stech} can be derived from our main framework.

We fix the notation as in Proposition~\ref{pro:stech-ps}. In particular, $\cal N$ is the set of all bounded seminorms on a Banach space $(E,\|\cdot\|)$. Obviously $\cal N$ is a convex cone: $\alpha\nu_1+\beta\nu_2\in \cal N$ for all $\nu_i\in \cal N$ and all $\alpha,\beta\ge0$.
The rigorous definition of the convex cone ${\cal N}_0\subset \cal N$ of equivalent norms on $E$ is:
$$
    \nu \in {\cal N}_0 \iff \exists\ a,b>0 \colon  a\|x\|\le\nu(x)\le b\|x\|,\ \forall x\in E.
$$
Since $\nu\in\cal N$ implies $\nu(x) \le \rho(\nu,0) \|x\|$ for all $x\in E$, we have
\begin{equation}
    \label{eq:a-nu}
    v\in {\cal N}_0 \iff a_\nu := \inf \nu (S_E) > 0.
\end{equation}
From this it is clear that
$$
    {\cal N} + {\cal N}_0 = {\cal N}_0,
$$
and, in particular, $\nu + \varepsilon\|\cdot\| \in {\cal N}_0$ for any $\nu\in\cal N$ and any $\varepsilon > 0$. So, ${\cal N}_0$ is dense in $\cal N$. ${\cal N}_0$ is also open, as it follows from the following lemma. Then, obviously, ${\cal N} = \overline {{\cal N}_0}$ and ${\cal N}_0 = {\cal N}^\circ$ (the interior of ${\cal N}$).
\begin{lem}
    \label{lem:N0-open}
    If $\nu\in {\cal N}_0$ then $B_{a_\nu}^\circ (\nu) \subset {\cal N}_0$.
\end{lem}
\begin{proof}
    Take a $\nu'\in {\cal N}$ such that $\rho(\nu',\nu) < a_\nu$. Let $h\in S_E$ be arbitrary. Then $\nu'(h) = \nu(h) + \nu'(h) - \nu(h) \ge a_\nu - \sup_{x\in B_X} |\nu'(x) - \nu(x)| = a_\nu - \rho(\nu',\nu)$. So, $a_{\nu'} \ge a_\nu - \rho(\nu',\nu) > 0$ and from \eqref{eq:a-nu} it follows that $\nu'\in {\cal N}_0$.
\end{proof}
\begin{lem}
    \label{lem:stech-key}
     Let $M\subset E$ be a non-empty closed, convex and bounded set. Let $\nu\in\cal N$ be such that
    $$
        \inf_M \nu  > 0.
    $$
    Then for each $\varepsilon>0$ there is $\nu'\in B_\varepsilon(\nu)$ and $\delta>0$ such that
    \begin{equation}
        \label{eq:stech-key}
        \diam \{x\in M:\ \nu'(x) \le \inf_M \nu'  + \delta\} < \varepsilon.
    \end{equation}
\end{lem}
\begin{proof}
    Fix $\varepsilon > 0$.

    Let $(x_n)_{n=1}^\infty\subset M$ be a minimizing sequence for $\nu$, i.e.
    $$
        \lim_{n\to\infty} \nu(x_n) = \inf_M \nu .
    $$
    Define $\nu_n\in\cal N$ by
    $$
        \nu_n(x) := \nu(x) + \varepsilon \min_{t\in\R}\|x-tx_n\|.
    $$
    Note that $\nu_n\in\cal N$, because $\tilde \nu_n(x)=\eps \min_{t\in\R}\|x-tx_n\|\in\cal N$. Also, since $\rho (0,\tilde \nu_n)\le\varepsilon$, we have $\nu_n\in B_\varepsilon(\nu)$.

    Fix a sequence $\delta_n \searrow 0$. We claim that
    \begin{equation}
        \label{eq:stech-key-seq}
        \lim_{n\to\infty} \diam \{x\in M:\ \nu_n(x)\le \inf_M \nu_n + \delta_n\} = 0.
    \end{equation}
    If this is true, we can take $\nu'=\nu_n$ for some large enough $n$.

    Assume that \eqref{eq:stech-key-seq} is false, that is, there is $c>0$ such that for infinitely many $n$'s the diameter of the set in the left hand side of \eqref{eq:stech-key-seq} is greater than $2c$. For avoiding subscripts and because dropping some $x_n$'s ruins nothing, let these diameters be greater than $  2c$ for all $n\in\N$. This means that for all $n\in\N$ there are
    \begin{equation}
        \label{eq:yn-cond}
        y_n\in M \text{ such that }\ \nu_n(y_n) \le \inf_M \nu_n + \delta_n,\text{ and } \|x_n-y_n\| > c.
    \end{equation}
    Let $t_n\in\R$ be such that $\|y_n-t_nx_n\| = \min_{t\in\R}\|y_n-tx_n\|$.

    We have that
    \begin{equation}\label{trd}
    \nu(y_n) \le \nu_n(y_n) \le \inf_M \nu_n + \delta_n \le \nu_n(x_n) + \delta_n = \nu (x_n) + \delta_n.
    \end{equation}
     Therefore,
    \begin{equation}
        \label{eq:yn-minimising}
        \lim_{n\to\infty} \nu(y_n) = \inf_M \nu.
    \end{equation}
    From \eqref{trd} we get that $\nu_n(y_n) \le \nu (x_n) + \delta_n$ and having in mind that $\nu_n(y_n) = \nu(y_n) + \varepsilon\|y_n-t_nx_n\|$, we obtain
    $$
        \varepsilon\|y_n-t_nx_n\| \le \nu (x_n) - \nu(y_n) + \delta_n\to0,\text{ as }n\to\infty.
    $$
    Hence,
    \begin{equation}
        \label{eq:yn-tn-0}
        \lim_{n\to\infty} \|y_n-t_nx_n\| = 0.
    \end{equation}
    Since $\nu$ is Lipschitz, \eqref{eq:yn-tn-0} imply $\nu(y_n)-\nu(t_nx_n)\to0$, as $n\to\infty$, that is, $\nu(y_n) - |t_n|\nu(x_n)\to0$, as $n\to\infty$. Because, see \eqref{eq:yn-minimising}, both $\nu(x_n)$ and $\nu(y_n)$ tend to $\inf_M\nu > 0$, we can divide and get $\nu(y_n)/\nu(x_n) - |t_n|\to0$, as $n\to\infty$, but $\nu(y_n)/\nu(x_n) \to 1$, so
    $$
        \lim_{n\to\infty} |t_n| = 1.
    $$
    If there were a subsequence $(t_{n_k})_{k=1}^\infty$ such that $t_{n_k}\to1$, as $k\to\infty$, then because of \eqref{eq:yn-tn-0}, we would get $\|y_{n_k}-x_{n_k}\|\to0$, as $k\to\infty$, which contradicts \eqref{eq:yn-cond}. Therefore,
    $$
        \lim_{n\to\infty} t_n = -1.
    $$
    Because $(x_n)_{n=1}^\infty$ is a  bounded sequence (being in $M$, for example) this implies $\|(y_n+x_n)-(y_n-t_nx_n)\| = |1+t_n|\|x_n\|\to0$, as $n\to\infty$. From \eqref{eq:yn-tn-0} we get $(x_n+y_n)\to0$ as $n\to\infty$. Since $M$ is convex, we have $(x_n+y_n)/2\in M$ and because $M$ is closed, we get $0\in M$, contradiction.
\end{proof}

\begin{proof}[\textbf{\emph{Proof of Proposition~\ref{pro:stech-ps}}}]
    We will check Definition~\ref{pert_sp}(i)-(iv). Fix $\nu\in\cal N$.

    For a fixed $p\in E$ the function $g_p(x) = \nu(p-x)$, defined for $x\in M$, is Lipschitz, because $\nu$ is, and bounded on the bounded set $M$. Thus,   Definition~\ref{pert_sp}(i) is immediate.

    Definition~\ref{pert_sp}(ii) follows from
    $$
        \|g_p-g_q\|_\infty = \sup_{x\in M} |\nu(p-x)-\nu(q-x)| \le \nu(p-q)\le k_\nu\|p-q\|.
    $$

    For Definition~\ref{pert_sp}(iii) note that
    $$\|g\|_\infty = \sup_{x\in M}\nu(p-x) \le \sup_{B_E} \nu .\sup_{x\in M}\|p-x\| = \sup_{x\in M}\|p-x\|\rho(g,0),$$
    so $c(p) = \sup_{x\in M}\|p-x\| < \infty$ does the job.

    Finally, we will verify  Definition~\ref{pert_sp}(iv). Fix $\varepsilon>0$. If $p\in M$ we can just take
    $$
        \nu' = \nu + \varepsilon\|\cdot\|
    $$
    and the corresponding $g'(\cdot)=\nu'(p-\cdot)$. Clearly, $\rho(g,g')=\varepsilon$ and $g'$ has a strong minimum at $p$, so we are done in this case.

    Now, let $p\not\in M$. We first add to $\nu$ a small multiple of the original norm $\|\cdot\|$ to get $\tilde \nu=\nu+\eps/2\|\cdot\|$ such that
    $$
        \inf_{p-M}\tilde\nu > 0.
    $$
    Applying  Lemma~\ref{lem:stech-key} to   $\tilde \nu$ and  the set $p-M$ for $\eps/2$ we get $\nu'\in B_{\eps/2}(\tilde \nu)\subset B_\varepsilon(\nu)$ and $\delta>0$ such that
    $$
        \diam\{x\in p-M:\ \nu'(x) \le \inf_{p-M} \nu' + \delta\} < \varepsilon/2,
    $$
    for this is what \eqref{eq:stech-key} translates to. But this is the same as
    $$
        \diam\{x\in M:\ \nu'(p-x) \le \inf_{x\in M} \nu'(p-x) + \delta\} < \varepsilon/2,
    $$
    hence,
    $$
        \diam\Omega_{g'}(\delta) < \varepsilon,
    $$
    where $g'=\nu'(p-\cdot)$, so $g'\in B_\varepsilon(g)$, and Definition~\ref{pert_sp}(iv) follows.
\end{proof}

\begin{proof}[\textbf{\emph{Proof of Theorem~\ref{thm:stech}}}] Let $|\cdot|$ be an equivalent norm on $E$. Fix $\varepsilon > 0$ such that $\varepsilon < a_{|\cdot|}$, see \eqref{eq:a-nu}.

    From Proposition~\ref{pro:stech-ps} and Theorem~\ref{thm:main} for $f_p\equiv0$ on $M$, there is a dense $G_\delta$ subset $G$ of $(\cal N,\rho)$ with the properties listed in Theorem~\ref{thm:main}.

    Set $\tilde g(p,x)=|p-x|$ and take $g\in G$ such that $\rho(\tilde g,g) < \varepsilon$, and let
    $$
        g(p-x) = \nu(p-x),\text{ where }\nu\in\cal N.
    $$
    Because of $a_{|\cdot|}>\eps>\rho(\tilde g,g)=\rho(|\cdot|,\nu)$,    from Lemma~\ref{lem:N0-open} it follows that $\nu$ is an equivalent norm, say $\nu = \n\cdot\n$.

    For $g\in G$, according to Theorem~\ref{thm:main}, there is a dense $G_\delta$ subset $V$ of $E$ such that $W\subset V$ and for all $p\in V$, $f_p+g_p$ attains its strong minimum at $x_p\in M$. But because $f_p\equiv0$ on $M$, we have
    $$
        g_p(x_p) = \min_M g_p
    $$
      in strong sense (that is, the minimization problem is well-posed), or
    $$
        \n p - x_p \n = \min_{x\in M} \n p - x \n = \dist _{\n\cdot\n} (p,M)
    $$
    in strong sense.

    The continuity of the $\n\cdot\n$-metric projection $p\to x_p$ on $p\in V$, follows from \eqref{eq:main-cont}.
\end{proof}

\section*{Acknowledgements}
We wish to express our sincere gratitude to Dr. Milen Ivanov  for his guidance  and for the invaluable help and constant support.

Thanks are due to Prof. Pando Georgiev for bringing our attention to his results in \cite{Pando-st}.


\begin{thebibliography}{99}

    \bibitem{beer}
    G. Beer, Topologies on closed and closed convex sets, Mathematics and its Applications \textbf{268}, Kluwer Academic Publishers, Dordrecht (1993).

    \bibitem{B-P}
    J. Borwein and D. Preiss, A smooth variational principle with applications to subdifferentiability of convex functions, Trans. Amer. Math. Soc. 303 (1987), 517--527.


    \bibitem{Cobzas}
    S. Cobza\c{s}, Geometric properties of Banach spaces and the existence of nearest and farthest points,
    Abstract and Applied Analysis, 2005  (2005), 3, 259--285.

    \bibitem{Danzig}
    G. B. Dantzig, J. Folkman, and N. Shapiro,
    On the continuity of the minimum set of a continuous function,
    Journal of Mathematical Analysis and Applications 17 (1967), 519--548.

    \bibitem{DGZ-article}
    R. Deville, G. Godefroy, and V. Zizler,
    A smooth variational principle with applications to Hamilton-Jacobi equations in infinite dimensions,
    Journal of Functional Analysis, 111 (1993), 1, 197--212.

    \bibitem{DGZ}
    R. Deville, G. Godefroy, and V. Zizler,  Smoothness and renormings in
    Banach spaces.  Pitman Monographs and Surveys in Pure and Appl.
    Math. \textbf{64}, Longman Scientific \& Technical, Harlow (1993).

    \bibitem{dev-pro} R. Deville and A. Proh\'azka,
    A parametric variational principle and residuality,
    Journal of Functional Analysis, 256 (2009), 11, 3568--3587.

    \bibitem{EfimovStechkin}
        N. Efimov and S. Ste\v{c}kin, Some properties of \v{C}eby\v{s}ev sets,
        Dokl. Akad. Nauk SSSR 118 (1958), 1, 17--19 (in Russian).

\bibitem{Pando-st}
 P. Georgiev, The strong Ekeland variational principle, the strong drop theorem and applications,
 Journal of Mathematical Analysis and Applications 131 (1988), 1--21.
 
    \bibitem{pando}
    P. Georgiev, Parametric Ekeland's variational principle,
    Applied Mathematics Letters 14 (2001), 6, 691--696.

    \bibitem{pando2}
    P. Georgiev, Parametric Borwein-Preiss variational principle and applications,
    Proceedings of the American Mathematical Society 133 (2005), 11, 3211--3225.
    

 
 
    \bibitem{iv-zla-jca}
    M. Ivanov and N. Zlateva, Perturbation method for variational problems, J. Conv. Anal. 19 (2012), 4, 1033--1042.

    \bibitem{iv-zla-jota}
    M. Ivanov and N. Zlateva, Perturbation method for non-convex integral
    functional, JOTA 157 (2013), 3, 737--748.

    \bibitem{rocwets}
    R. T. Rockafellar and R. J.-B. Wets, Variational analysis, Grundlehren
    der mathematischen Wissenschaften, vol. 317. Springer, New York,
    1998.

 \bibitem{Stechkin}
   S. B. Ste\v{c}hkin,  Approximative properties of sets in normed  linear spaces, Revue Math. Pures Appl. 8 (1963), 5--18 (in Russian).


    \bibitem{Tbook}
    L. Thibault, Unilateral Variational Analysis in Banach Spaces.
    Part I: General Theory.
    Part II: Special Classes of Functions and Sets,  World Scientific, 2023, ISBN: 978-981-125-816-9.

    \bibitem{pertOrl}
    H. Topalova and N. Zlateva, Perturbation Method in Orlicz Sequence Spaces, Set-Valued and Variational Analysis, 32 (2024),  2, art. 12.

    \bibitem{SimultOptim}
    H. Topalova and N. Zlateva, Simultaneous perturbed minimization of a convergent sequence of functions. Optimization, (2024), 1--11, https://doi.org/10.1080/02331934.2024.2386112.

    \bibitem{vesely}
    L. Vesel\'y, A parametric smooth variational principle and support properties of convex sets and functions,
    Journal of Mathematical Analysis and Applications 350 (2019), 2, 550--561.

    \end{thebibliography}
\end{document}